\documentclass[journal,twocolumn]{IEEEtran}

\usepackage[english]{babel}
\usepackage{amsmath,amssymb,mathtools}
\let\labelindent\relax
\usepackage{enumitem}
\usepackage{pgfplotstable}
\usetikzlibrary{arrows.meta}
\usepackage{tikz}
\usepackage[ruled,vlined]{algorithm2e}
\SetKwBlock{FEND}{for $k=0,1,2,\ldots$ each node $i$ does}{}
\newcommand{\ubar}[1]{\underline{#1}}

\usepackage{booktabs}

\newenvironment{Protocol}[1][htb]
 {
 \begin{algorithm}[#1]%
 }{\end{algorithm}}

\newtheorem{definition}{Definition}
\newtheorem{problem}{Problem}

\newtheorem{theorem}{Theorem}
\newtheorem{corollary}{Corollary}

\newtheorem{assumption}{Assumption}

\newcommand{\rea}{\mathbb{R}}

\ifCLASSINFOpdf
\else
\fi
\begin{document}

\title{Dynamic Max-Consensus and Size Estimation of Anonymous Multi-Agent Networks}

\author{Diego Deplano, Mauro Franceschelli, Alessandro Giua 
\thanks{D. Deplano, M. Franceschelli and A. Giua are with DIEE, University of Cagliari, 09123 Cagliari, Italy.
Emails: {\tt \{diego.deplano,mauro.franceschelli,giua\}@unica.it}}
\thanks{Mauro Franceschelli is corresponding author.
This work was supported in part by the Italian Ministry of Research and Education (MIUR) with the grant "CoNetDomeSys", code RBSI14OF6H, under call SIR 2014 and by Region Sardinia (RAS) with project MOSIMA, RASSR05871, FSC 2014-2020, Annualità 2017, Area Tematica 3, Linea d'Azione 3.1.}
}

\maketitle

\begin{abstract}

In this paper we propose a novel consensus protocol for discrete-time multi-agent systems (MAS), which solves the dynamic consensus problem on the max value, i.e., the dynamic max-consensus problem. In the dynamic max-consensus problem to each agent is fed a an exogenous reference signal, the objective of each agent is to estimate the instantaneous and time-varying value of the maximum among the signals fed to the network, by exploiting only local and anonymous interactions among the agents. 
The absolute and relative tracking error of the proposed distributed control protocol is theoretically characterized and is shown to be bounded and by tuning its parameters it is possible to trade-off convergence time for steady-state error.
The dynamic Max-consensus algorithm is then applied to solve the distributed size estimation problem in a dynamic setting where the size of the network is time-varying during the execution of the estimation algorithm. Numerical simulations are provided to corroborate the theoretical analysis. 
\end{abstract}


%
\IEEEpeerreviewmaketitle

\section{Introduction}

Recent years have seen considerable interest in the design of algorithms to solve the problem of consensus agreement over networked systems. The problem is considered to be solved when the agents agree upon a value while using only local information about neighboring agents. Most of the literature usually considers static agreement values, i.e., agents are required to converge to a static reference signal. On the contrary, in the dynamic consensus problem a time-varying reference signal is associated to each agent, and they are required to achieve consensus upon a function of the time-varying reference signals, such as average, median, maximum and so on. 


While the literature has focused significantly on the dynamic average-consensus problem \cite{spanos2005dynamic,Freeman2006,Zhu2010322,Montijano20143131,FRANCESCHELLI201969,Freeman2010,Freeman2015,Cortes2015,FranOpen1}, estimating the average is not the only attractive goal. In particular, the focus of this work is the the development of dynamic protocols achieving consensus upon on the max value among the reference signals. Applications of dynamic max-consensus protocols mainly reside in the field of distributed synchronization, such as time-synchronization \cite{Dengchang} and target tracking \cite{Petitti11}, and network parameter estimation, such as cardinality \cite{Lucchese15} and highest/lowest node degree \cite{Borsche10}.

\textbf{Literature review.} In the literature the so-called max-consensus problem has been thoroughly investigated. Its objective is to make the states of the agents converge to the maximum of their initial states. 
The most popular max-consensus protocol consists in initializing the network to a set of values and let agent update its state at each instant of time by taking the maximum value among the value of the neighbors' state and its own state \cite{Olfati04}. 
The work in \cite{Nejad09} proposes conditions to achieve max-consensus and compute convergence rate of these protocols for different communication topologies. 
Only little effort has been paid to analyze slightly different but much more complicated variations of this problem. In particular, convergence results have only been provided for synchronous switching topologies \cite{Nejad10} and for probabilistic
asynchronous fixed frameworks \cite{Iutzeler12}. The contribution of introducing time delays in the communications is due to \cite{Giannini16}, while \cite{Zhang16} is the first work allowing noise in the communications. Finally, the case with agents with the possibility to join or leave the network, so-called open multi-agent systems, is addressed in \cite{Abdelrahim}.
If a static consensus protocol is used to perform a distributed estimation upon some time-varying quantities, known or measured by the agents, the protocol requires to be re-initialized in the whole network each time the value of the function to be estimated changes. Dynamic consensus protocols has been introduced in order to overcome the issue of re-initializing the network.

\textbf{Main contribution.} Up to the authors' knowledge, there no exists dynamic max-consensus protocol in the literature and, moreover, the existing results to solve the max-consensus problem only apply to static networks. Therefore, in this work we propose two protocols capable of solving the dynamic max-consensus problem in time-varying networks which are connected at any time. Furthermore, variants of these protocols are proposed to solve the dual problem of estimating the min value.

\section{Open Dynamical Systems}\label{sec:preliminaries}

We consider Multi-Agent Systems (MAS) with a time-varying number of agents. These MAS are also called \emph{Open} MAS (OMAS). At any time $k\in\mathbb{N}$, let $\mathcal{G}_k=(V_k,E_k) $ be the undirected graph describing the pattern of interactions among agents, where $V_k\subset\mathbb{N}$ is the set of active nodes at time $k$ and $E_k\subseteq (V_k\times V_k)$ is the set of active communication channels between agents in the network. The number $n_k$ of active agents in the network at time $k$ is given by the cardinality of the set $V_k$, i.e., $n_k=|V_k|$. Agents $p$ and $q$ are said to be neighbors at time $k$ if there exists an edge from $p$ to $q$ (and vice versa) at time $k$, i.e., $(p,q),(p,q)\in E_k$. A set of neighbors $\mathcal{N}^p_k$ is associated to each node $p$ at time $k$, defined as $\mathcal{N}^p_k=\left\{q\in V_k: (p,q)\in E_k\right\}$, which represents the agents in the graph which share a point-to-point communication channel with agent $p$.


A \emph{path} $\pi^{pq}_k$ between two nodes $p$ and $q$ in a graph is a finite sequence of $m$ edges $e_\ell=(i_\ell,j_\ell)\in E_k$ that joins node $p$ to node $q$, i.e., $i_1=p$, $j_m=q$ and $j_\ell=i_{\ell+1}$ for $\ell=1,\ldots,m-1$. An undirected graph is said to be \emph{connected} at time $k$ if there exists a path $\pi^{pq}_k$ between any pair of nodes $p,q \in V_k$. The \emph{diameter} of a graph, denoted as $\delta(\mathcal{G}_k)$, is defined as the longest among the shortest paths among any pair of nodes $p,q\in V_k$. For any connected undirected graph it holds $\delta(\mathcal{G}_k)\leq n_k-1$.

For each time $k$ and agent $p\in N_k$, the scalar agent's state is denoted $x_k^p\in\mathbb{R}$ and its input is denoted $u_k^p$, which are defined only at time instants such that $p\in V_k$. More generally, in this paper we shall call \emph{open sequence} any sequence $\{y_k:k\in\mathbb{N}\}$, or simply $\{y_k\}$, where $y_k\in\mathbb{R}^{V_k}$. When the number of nodes $n_k$ changes with the time $k$, it is not possible in general to write $x_{k+1}$ as a function solely of $x_k$: therefore, the evolution of $x_k$ does not constitute a "closed" dynamical system. Thus, one needs to first partition the nodes into three sets:
\begin{itemize}
\item The \emph{departing} nodes $D_k = V_k \setminus V_{k+1}$;
\item The \emph{arriving} nodes $A_k = V_{k+1} \setminus V_{k}$;
\item The \emph{remaining} nodes $R_k = V_k \cap V_{k+1}$.
\end{itemize}
With this notation, we shall define the evolution of each agent depending on its status,
\begin{equation}\label{eq:open_loc_mas_var}
x_{k+1}^p = 
\begin{cases}
\emptyset & \text{if } p\in D_k\\
f_a(x_{k+1},V_{k+1},E_{k+1},u_{k+1}) & \text{if } p\in A_k\\
f_r(x_k,V_k,E_k,u_k) & \text{if } p\in R_k
\end{cases}.
\end{equation}
In other words, since $x_{k+1}$ must take values in $\mathbb{R}^{V_{k+1}}$ for all $k\in\mathbb{N}$, the departing nodes have to be left out, the arriving nodes need to be initialized when joining to the network according to some rule $f_a(\cdot)$ and the remaining nodes update their current state according to some rule $f_r(\cdot)$.
We observe that if the set of agents does not change at time $k+1$, i.e., $V_k=V_{k+1}$, then we can write in vector form the dynamics of the corresponding MAS as
\begin{equation}\label{eq:open_glob_mas_stat}
x_{k+1} = f(x_k,V_k,E_k,u_k).
\end{equation}
\section{Problem statement} \label{sec:problemstatement}

The dynamic consensus problem in a MAS as in \eqref{eq:open_glob_mas_stat} consists in steering the agents' states to track a function $ g(\cdot) : \mathbb{R}^n\rightarrow \mathbb{R}$ of the time-varying exogenous reference signals $ u^i_k $ such that the tracking error
\begin{equation}\label{eq:trackingerror}
e_k=\max_{i\in V_k} |x^i_k-g(u_k)|
\end{equation}
converges within a boundary layer for any initial condition. The \emph{convergence time} $T^c$ denotes the time required by the MAS to achieve such a bounded tracking error.

If the MAS is open, i.e., agents may join or leave the network, and the natural way to generalize the dynamic consensus problem to OMAS is requiring that at every change of the network the tracking error \eqref{eq:trackingerror} decrease before another change, eventually converging within a boundary layer. Since the change in the network is a local event, it is not possible to re-initialize the whole network in a distributed way, i.e., the protocol must be robust to the initial condition. Furthermore, the change may results in a unpredictable discontinuity of the quantity of interests, thus requiring the whole network to agree upon a value which is uncorrelated to the one at the previous step. In the light of the above considerations, we have to make some assumptions on the rate of change of the network and the inputs.

We consider a network of agents whose topology at time $k\in\mathbb{N}$ is represented by an undirected connected graph ${\mathcal{G}_k=(V_k,E_k)}$ satisfying the next assumption.
{\begin{assumption}\label{ass:slow_net}
There exists an interval time $\Upsilon\in\mathbb{N}$ such that two changes of the network $\mathcal{G}$ can not occur in $\Upsilon$ units of time, i.e., 
\begin{align*}\label{eq:slow_net}
\forall k_0\in\mathbb{N}:
\:\mathcal{G}_{k_0-1} \neq \mathcal{G}_{k^*}
\: \Rightarrow \:
\mathcal{G}_k = \mathcal{G}_{k+1},
\:\forall k \in\{k_0,k_0+\Upsilon\}.
\end{align*}
\hfill $\blacksquare$
\end{assumption}\medskip}

Each agent $i$ has access to a time-varying external reference signal $u^i_k\in\mathbb{R}$ satisfying the next assumption.
{\begin{assumption}\label{ass:boundedu}
Each unknown exogenous reference signal is such that its has absolute change is bounded by a constant $\Pi\in \mathbb{R}_+$, i.e., 
\begin{equation}\label{eq:boundedu}
|u^i_{k+1}-u^i_k| \leq \Pi, \quad \forall i \in V_k, \; \forall k\in \mathbb{N}.
\end{equation}
\hfill $\blacksquare$
\end{assumption}\medskip}
Note that, due to Assumptions~\ref{ass:slow_net}-\ref{ass:boundedu}, within two changes, the OMAS can be regarded as a MAS with dynamics given in \eqref{eq:open_glob_mas_stat}, but the inputs keep varying during such an interval. In such an interval one can define the \emph{transient time} and the \emph{convergence time} as follows.

\begin{definition}\label{def:trantime}
The \emph{transient time} $T^t_{k_0}\in\mathbb{N}$ is the time that the network needs to remain unchanged in order achieve a decreasing tracking error \eqref{eq:trackingerror} from an initial time $k_0\in\mathbb{N}$.
\end{definition}
\begin{definition}\label{def:convtime}
The \emph{convergence time} $T^c_{k_0}\in\mathbb{N}$ is the time that the network needs to remain unchanged in order achieve a bounded tracking error \eqref{eq:trackingerror} from an initial time $k_0\in\mathbb{N}$.\hfill $\blacksquare$
\end{definition}

In particular, in this work we focus on the so-called dynamic min/max-consensus problem, i.e., the quantity of interest is either the minimum $\ubar{u}_k$ or the maximum $\bar{u}_k$ of the the exogenous reference signals,  defined as follows
\begin{equation}\label{eq:ustar}
\ubar{u}_k=\min_{i\in V_k} u^i_k,\qquad \bar{u}_k=\max_{i\in V_k} u^i_k.
\end{equation}

\begin{problem}\label{prb:main}
Consider an OMAS as in \eqref{eq:open_loc_mas_var} satisfying Assumptions~\ref{ass:slow_net}-\ref{ass:boundedu} and let the quantity of interest $g(\cdot)$ be either the maximum $\bar{u}$ or the minimum $\ubar{u}$ as in \eqref{eq:ustar}.

The dynamic min/max-consensus problem consists in designing the local interaction rules $f_a(\cdot)$, $f_r(\cdot)$ such that the tracking error \eqref{eq:trackingerror} satisfies
\begin{align}
e_{k+1}&< e_k & \text{for } k\in [k_0+T^t_{k_0},k_0+\min\{T^c_{k_0},\Upsilon\}],\label{eq:err_dec}\\
e_{k+1}&\leq \varepsilon & \text{for } k\in [k_0+\min\{T^c_{k_0},\Upsilon\},k_0+\Upsilon],\label{eq:err_lay}
\end{align}
for any $k_0\in\mathbb{N}$ such that $\mathcal{G}_{k_0-1}\neq \mathcal{G}_{k_0}$ and where
\begin{itemize}
\item $T^t_{k_0}\in[0,\min\{T^c_{k_0},\Upsilon\})$ is the transient time;
\item $T^c_{k_0}\in\mathbb{N}$ is the transient time;
\item $\varepsilon\in\mathbb{R}_+$ is the bound on the tracking error.
\end{itemize}
\hfill $\blacksquare$
\end{problem}
Problem~\ref{prb:main} requires two ingredients for declaring the dynamic min/max-consensus problem solved, which are the following:
\begin{itemize}
\item \emph{Condition \eqref{eq:err_dec}}: The tracking error must decrease between two changes in the network after a transient time.
\item \emph{Condition \eqref{eq:err_lay}}: The tracking error is bounded if the rate of change of the network is large enough, i.e., $\Upsilon\geq T^c_{k_0}$.
\end{itemize}

Objective of this paper is thus to propose local interaction protocols~\eqref{eq:open_loc_mas_var} for a discrete-time OMAS, which solve the dynamic consensus problem formalized in Problem~\ref{prb:main} for anonymous networks, i.e., the agents do not share their identification. The proposed protocols are then applied to solve the distributed online estimation problem of size the time-varying network to which the agents belong.

\section{Dynamic Max-Consensus Protocols} \label{sec:DMC}

In this section we provide two protocols to solve the dynamic max-consenus problem:
\begin{enumerate}
\item \emph{Approximate Dynamic Max-Consensus} (ADMC) Protocol: it enables the agents to converge to an \emph{approximate} consensus on the max value without requiring any information about the network graph.
\item \emph{Exact Dynamic Max-Consensus} (EDMC) Protocol: it enables the agents to reach an \emph{exact} consensus on the max value by requiring the knowledge of an upperbound of the network's graph diameter.
\end{enumerate}

For each protocol we prove that Problem~\ref{prb:main} is solved for $g(u_k)=\bar{u}_k$ as in \eqref{eq:ustar}, characterizing the convergence time $T^c_{k_0}$ as in Definition~\ref{def:convtime} and the bound on the tracking error
\begin{equation}\label{eq:abserr_max}
e(k)=\max_{i\in V} |x_i(k)-\bar{u}(k)|.
\end{equation}

\subsection{Approximate Consensus}\label{sec:DMC:A}

\begin{Protocol}[t]
 \SetKwInOut{Input}{Input}
 \SetKwInOut{Output}{Output}
 \SetKwInOut{Set}{Set}
 \Input{Tuning parameter $\alpha \in (0,1)$.\\
 }
\Output{Current state $x^i_k\in\mathbb{R}$ for $i\in V_k$.}
\FEND{
\If{$i\in A_k$\textbf{ \emph{or}} $k=0$}{
$x^i_{k+1}\gets u^i_{k+1}$
 }
\ElseIf{$i\in R_k$}{
Gather $x^j_k$ from each neighbor $j \in \mathcal{N}^i_k$\\
Update the current state according to \\
$\displaystyle x^i_{k+1} \gets \max_{j \in \mathcal{N}^i_k\bigcup \left\{i\right\}} \left\{x^j_k - \alpha,u^i_k\right\}$\\
}
}
\caption{{\small Approximate Dynamic Max-Consensus (ADMC)}}
\label{alg:ADMC}
\end{Protocol}

In Protocol~\ref{alg:ADMC} we detailed the ADMC Protocol, which makes use of the following local interaction rule:

\begin{equation}\label{eq:protocol1}
x^i_{k+1} =
\begin{cases}
\displaystyle\max_{j \in \mathcal{N}^i_k\bigcup \left\{i\right\}} \left\{ x^j_k-\alpha,u^i_k\right\}& \text{if } i\in R_k\\
u^i_{k+1} & \text{if } i\in A_k
\end{cases}
\end{equation}

where $\alpha\in(0,1)$ is a scalar tuning parameter.  At each iteration all remaining agents gather the state values of their neighbors and update their state according to~\eqref{eq:protocol1}, which only requires local communications, while all arriving agents initialize their next state at their own inputs. Note that at the initial time $k=0$ all agents are arriving agents.

\begin{theorem}[ADMC Protocol: Tracking Error]\label{th:ADMC_track}
Consider an OMAS executing Protocol~\ref{alg:ADMC} with tuning parameter $\alpha \in \rea_+$ under Assumptions~\ref{ass:slow_net}-\ref{ass:boundedu}. Consider a generic time $k_0\in\mathbb{N}$ at which the network changes, i.e., $\mathcal{G}_{k_0-1}\neq \mathcal{G}_{k_0}$.

If $\mathcal{G}_{k_0}$ is connected with diameter $\delta(\mathcal{G}_{k_0})$ and if 
\begin{equation}\label{eq:ADMC_alpha}
\alpha> \Pi,\quad \Upsilon\geq \delta(\mathcal{G}_{k_0})
\end{equation} 
Problem~\ref{prb:main} is solved with transient and convergence times
\begin{align}\label{eq:ADMC_T}
T^t_{k_0} &= \delta(\mathcal{G}_{k_0}),\\
T^c_{k_0} &= \max \left\{ T^t_{k_0},\left\lceil \frac{\max \left\{x_{k_0}-\bar{u}_{k_0}\right\}}{\alpha - \Pi}\right\rceil\right\},\nonumber
\end{align}
with $\bar{u}_k$ defined in \eqref{eq:ustar}. Moreover, if $\Upsilon \geq T^c_{k_0}$, the tracking error \eqref{eq:abserr_max} is bounded for any $k\in[k_0+T^c_{k_0},k_0+\Upsilon]$ by
\begin{equation}\label{eq:ADMC_err}
e_k\leq \varepsilon = (\delta(\mathcal{G}_{k_0})+1)\Pi+\alpha \delta(\mathcal{G}_{k_0}).
\end{equation}
\end{theorem} \medskip

\begin{IEEEproof}
The proof is given in Appendix.
\end{IEEEproof}

\begin{corollary}[ADMC Protocol: Steady State Error]\label{th:ADMC_steady}
Consider an OMAS executing Protocol~\ref{alg:ADMC} with tuning parameter $\alpha \in \rea_+$ under Assumption~\ref{ass:slow_net}. Consider a generic time $k_0\in\mathbb{N}$ at which the network changes, i.e., $\mathcal{G}_{k_0-1}\neq \mathcal{G}_{k_0}$ and let the inputs be constant for $k\in[k_0,k_0+\Upsilon]$.

If $\mathcal{G}_{k_0}$ is connected with diameter $\delta(\mathcal{G}_{k_0})$ and if 
\begin{equation*}
\alpha> \Pi,\quad \Upsilon\geq \delta(\mathcal{G}_{k_0})
\end{equation*} 
Problem~\ref{prb:main} is solved with transient and convergence times
\begin{align}\label{eq:ADMC_T}
T^t_{k_0} &= \delta(\mathcal{G}_{k_0}),\\
T^c_{k_0} &= \max \left\{ T^t_{k_0},\left\lceil \frac{\max \left\{x_{k_0}-\bar{u}_{k_0}\right\}}{\alpha}\right\rceil\right\},\nonumber
\end{align}
with $\bar{u}(k)$ defined in \eqref{eq:ustar}. Moreover, if $\Upsilon\geq T^c_{k_0}$, the tracking error \eqref{eq:abserr_max} is fixed for any $k\in[k_0+T^c_{k_0},k_0+\Upsilon]$
\begin{equation}\label{eq:ADMC_err_ss}
e_k \leq \varepsilon_{ss} = \alpha \delta(\mathcal{G}_{k_0})).
\end{equation}
\end{corollary} \medskip
\begin{IEEEproof} Since the inputs are constant for $k\in[k_0,k_0+\Upsilon]$, Assumption~\ref{ass:boundedu} is satisfied with $\Pi=0$. Therefore, we can apply Theorem~\ref{th:ADMC_track} and specialize transient and convergence times given in \eqref{eq:ADMC_T} along with the tracking error given in \eqref{eq:ADMC_err} for $\Pi=0$, completing the proof of the corollary. Furthermore the bound is now a strict condition. In fact, since the inputs are constant, the network reaches an equilibrium and so the steady state error does not change over time.
\end{IEEEproof}

From the result of Theorem~\ref{th:ADMC_track} it follows that, according to~\eqref{eq:ADMC_err}, to minimize the absolute estimation error we need to choose $\alpha\approx 0$, $\alpha>\Pi\geq 0$. On the other hand, $\alpha$ determines the convergence time $T$ according to~\eqref{eq:ADMC_T}, with smaller values of $\alpha$ giving a greater convergence time. Thus, the value of $\alpha$ trades-off estimation error and convergence time. 

It follows that a pragmatic design criterion for the choice of $\alpha$ is to first fix the desired steady-state error and then choose the largest $\alpha$ which allows to satisfy the error performance constraint while minimizing the convergence time.

\subsection{Exact Consensus}\label{sec:DMC:E}

\begin{Protocol}[t]
 \SetKwInOut{Input}{Input}
 \SetKwInOut{Output}{Output}
 \SetKwInOut{Set}{Set}
 \Input{Network's diameter upperbound $\Delta \in \mathbb{N}$.}
\Output{Current state $x_{i\Delta}(k)\in\mathbb{R}$ for $i\in V$.}
\FEND{
\If{$i\in A_k$\textbf{ \emph{or}} $k=0$}{
\For{$\ell =0,1,\ldots,\Delta$}{
$x^{i\ell}_{k+1}\gets u^{i\ell}_{k+1}$
}
 }
\ElseIf{$i\in R_k$}{
Gather $[x^{j0}_k,\ldots,x^{j\Delta}_k]$ from each neighbor $j \in \mathcal{N}^i$\\
Update the current states according to \\
$\displaystyle x^{i0}_{k+1} \gets u^i_k$\\
\For{$\ell =1,\ldots,\Delta$}{
$\displaystyle x^{i\ell}_{k+1} \gets \max_{j \in \mathcal{N}^i_k\bigcup \left\{i\right\}} x^{j(\ell-1)}_k$
}}}
\caption{{\small Exact Dynamic Max-Consensus (EDMC)}}
\label{alg:EDMC}
\end{Protocol}

In Protocol~\ref{alg:EDMC} we detailed the EDMC Protocol, which makes use of the following local interaction rule
\begin{align}\label{eq:protocol2}
\displaystyle x^{i0}_{k+1} &= 
\begin{cases}
u^i_{k}& \text{if } i\in R_k\\
u^i_{k+1} & \text{if } i\in A_k
\end{cases}
\\
x^{i\ell}_{k+1} &=
\begin{cases}
\displaystyle \max_{j \in \mathcal{N}^i_k\bigcup \left\{i\right\}} x^{j(\ell-1)}_k& \text{if } i\in R_k\\
u^i_{k+1} & \text{if } i\in A_k
\end{cases}\nonumber
\end{align}
with $\ell=1,\ldots,\Delta$ and $\Delta\in\mathbb{N}$ is an upperbound on the the diameter of the underlying comunication network, i.e., $\delta \geq \delta(\mathcal{G}_k)$ at any time $k\in \mathbb{N}$. At each iteration all agents gather the state values of their neighbors and update their state according to~\eqref{eq:protocol2}, which only requires local communications, while all arriving agents initialize their next state at their own inputs. Note that at the initial time $k=0$ all agents are arriving agents.

\begin{theorem}[EDMC Protocol: Tracking Error]\label{th:EDMC_track}
Consider an OMAS executing Protocol~\ref{alg:EDMC} under Assumptions~\ref{ass:slow_net}-\ref{ass:boundedu}. Consider a generic time $k_0\in\mathbb{N}$ at which the network changes, i.e., $\mathcal{G}_{k_0-1}\neq \mathcal{G}_{k_0}$.

If graph $\mathcal{G}_{k_0}$ is connected with diameter $\delta(\mathcal{G}_{k_0})$ and if 
\begin{equation}\label{eq:EDMC_delta}
\Upsilon \geq \Delta \geq \delta(\mathcal{G}_{k_0}),
\end{equation} 
Problem~\ref{prb:main} is solved with transient and convergence times
\begin{equation}\label{eq:EDMC_T}
T^c_{k_0} = T^t_{k_0} = \Delta.
\end{equation}
Moreover, if $\Upsilon \geq T^c_{k_0}$, the tracking error \eqref{eq:abserr_max} is bounded for any $k\in[k_0+T^c_{k_0},k_0+\Upsilon]$ by
\begin{equation}\label{eq:EDMC_err}
e_k \leq \varepsilon = (\Delta + 1) \Pi.
\end{equation}
\end{theorem} \medskip
\begin{IEEEproof}
The proof is given in Appendix.
\end{IEEEproof}

\begin{corollary}[EDMC Protocol: Steady State Error]\label{th:EDMC_steady}
Consider an OMAS executing Protocol~\ref{alg:EDMC} under Assumption~\ref{ass:slow_net}. Consider a generic time $k_0\in\mathbb{N}$ at which the network changes, i.e., $\mathcal{G}_{k_0-1}\neq \mathcal{G}_{k_0}$ and let the inputs be constant for $k\in[k_0,k_0+\Upsilon]$.

If graph $\mathcal{G}_{k_0}$ is connected with diameter $\delta(\mathcal{G}_{k_0})$ and if 
\begin{equation*}
\Upsilon \geq \Delta \geq \delta(\mathcal{G}_{k_0}),
\end{equation*} 
Problem~\ref{prb:main} is solved with transient and convergence times
\begin{equation}\label{eq:EDMC_T_ss}
T^c_{k_0} = T^t_{k_0} = \Delta.
\end{equation}
Moreover, if $\Upsilon \geq T^c_{k_0}$, the tracking error \eqref{eq:abserr_max} is bounded for any $k\in[k_0+T^c_{k_0},k_0+\Upsilon]$ by
\begin{equation}\label{eq:EDMC_err_ss}
e_k = \varepsilon_{ss} = 0.
\end{equation}
\end{corollary} \medskip

\begin{IEEEproof} The proof is similar to the one of Corollary \ref{th:ADMC_steady}.
\end{IEEEproof}


\section{Dynamic Min-Consensus Protocols} \label{sec:DmC}

In this section we provide two protocols to solve the dynamic min-consenus problem:
\begin{enumerate}
\item \emph{Approximate Dynamic Min-Consensus} (ADmC) Protocol: it enables the agents to converge to an \emph{approximate} consensus on the min value without requiring any information about the network graph.
\item \emph{Exact Dynamic Min-Consensus} (EDmC) Protocol: it enables the agents to reach an \emph{exact} consensus on the min value by requiring the knowledge of an upperbound of the network's graph diameter.
\end{enumerate}

For each protocol we prove that Problem~\ref{prb:main} is solved for $g(u_k)=\ubar{u}_k$ as in \eqref{eq:ustar}, characterizing the convergence time $T^c_{k_0}$ as in Definition~\ref{def:convtime} and the bound on the tracking error
\begin{equation}\label{eq:abserr_min}
e(k)=\max_{i\in V} |x_i(k)-\ubar{u}(k)|.
\end{equation}

\begin{Protocol}[t]
 \SetKwInOut{Input}{Input}
 \SetKwInOut{Output}{Output}
 \SetKwInOut{Set}{Set}
 \Input{Tuning parameter $\alpha \in (0,1)$.\\
 }
\Output{Current state $x^i_k\in\mathbb{R}$ for $i\in V_k$.}
\FEND{
\If{$i\in A_k$\textbf{ \emph{or}} $k=0$}{
$x^i_{k+1}\gets u^i_{k+1}$
 }
\ElseIf{$i\in R_k$}{
Gather $x^j_k$ from each neighbor $j \in \mathcal{N}^i_k$\\
Update the current state according to \\
$\displaystyle x^i_{k+1} \gets \min_{j \in \mathcal{N}^i_k\bigcup \left\{i\right\}} \left\{x^j_k + \alpha,u^i_k\right\}$\\
}
}
\caption{{\small Approximate Dynamic Min-Consensus (ADmC)}}
\label{alg:ADmC}
\end{Protocol}

\subsection{Approximate Consensus}\label{sec:DmC:A}

In Protocol~\ref{alg:ADmC} we detailed the ADmC Protocol, which makes use of the following local interaction rule:

\begin{equation}\label{eq:protocol1}
x^i_{k+1} =
\begin{cases}
\displaystyle\min_{j \in \mathcal{N}^i_k\bigcup \left\{i\right\}} \left\{ x^j_k+\alpha,u^i_k\right\}& \text{if } i\in R_k\\
u^i_{k+1} & \text{if } i\in A_k
\end{cases}
\end{equation}

where $\alpha\in(0,1)$ is a scalar tuning parameter.  At each iteration all remaining agents gather the state values of their neighbors and update their state according to~\eqref{eq:protocol1}, which only requires local communications, while all arriving agents initialize their next state at their own inputs. Note that at the initial time $k=0$ all agents are arriving agents.

\begin{theorem}[ADMC Protocol: Tracking Error]\label{th:ADmC_track}
Consider an OMAS executing Protocol~\ref{alg:ADmC} with tuning parameter $\alpha \in \rea_+$ under Assumptions~\ref{ass:slow_net}-\ref{ass:boundedu}. Consider a generic time $k_0\in\mathbb{N}$ at which the network changes, i.e., $\mathcal{G}_{k_0-1}\neq \mathcal{G}_{k_0}$.

If $\mathcal{G}_{k_0}$ is connected with diameter $\delta(\mathcal{G}_{k_0})$ and if 
\begin{equation}\label{eq:ADmC_alpha}
\alpha> \Pi,\quad \Upsilon\geq \delta(\mathcal{G}_{k_0})
\end{equation} 
Problem~\ref{prb:main} is solved with transient and convergence times
\begin{align}\label{eq:ADmC_T}
T^t_{k_0} &= \delta(\mathcal{G}_{k_0}),\\
T^c_{k_0} &= \max \left\{ T^t_{k_0},\left\lceil \frac{\max \left\{\ubar{u}_{k_0}-x_{k_0}\right\}}{\alpha - \Pi}\right\rceil\right\},\nonumber
\end{align}
with $\bar{u}_k$ defined in \eqref{eq:ustar}. Moreover, if $\Upsilon \geq T^c_{k_0}$, the tracking error \eqref{eq:abserr_max} is bounded for any $k\in[k_0+T^c_{k_0},k_0+\Upsilon]$ by
\begin{equation}\label{eq:ADmC_err}
e_k\leq \varepsilon = (\delta(\mathcal{G}_{k_0})+1)\Pi+\alpha \delta(\mathcal{G}_{k_0}).
\end{equation}
\end{theorem} \medskip
 \begin{IEEEproof}
The proof is given in Appendix.
\end{IEEEproof}

\begin{corollary}[ADMC Protocol: Steady State Error]\label{th:ADmC_steady}
Consider an OMAS executing Protocol~\ref{alg:ADmC} with tuning parameter $\alpha \in \rea_+$ under Assumption~\ref{ass:slow_net}. Consider a generic time $k_0\in\mathbb{N}$ at which the network changes, i.e., $\mathcal{G}_{k_0-1}\neq \mathcal{G}_{k_0}$ and let the inputs be constant for $k\in[k_0,k_0+\Upsilon]$.

If $\mathcal{G}_{k_0}$ is connected with diameter $\delta(\mathcal{G}_{k_0})$ and if 
\begin{equation*}
\alpha> \Pi,\quad \Upsilon\geq \delta(\mathcal{G}_{k_0})
\end{equation*} 
Problem~\ref{prb:main} is solved with transient and convergence times
\begin{align}\label{eq:ADmC_T}
T^t_{k_0} &= \delta(\mathcal{G}_{k_0}),\\
T^c_{k_0} &= \max \left\{ T^t_{k_0},\left\lceil \frac{\max \left\{\ubar{u}_{k_0}-x_{k_0}\right\}}{\alpha}\right\rceil\right\},\nonumber
\end{align}
with $\bar{u}(k)$ defined in \eqref{eq:ustar}. Moreover, if $\Upsilon\geq T^c_{k_0}$, the tracking error \eqref{eq:abserr_max} is bounded for any $k\in[k_0+T^c_{k_0},k_0+\Upsilon]$ by
\begin{equation}\label{eq:ADmC_err_ss}
e_k \leq \varepsilon_{ss} = \alpha \delta(\mathcal{G}_{k_0})).
\end{equation}
\end{corollary} \medskip
\begin{IEEEproof} Since the inputs are constant for $k\in[k_0,k_0+\Upsilon]$, Assumption~\ref{ass:boundedu} is satisfied with $\Pi=0$. Therefore, we can apply Theorem~\ref{th:ADmC_track} and specialize transient and convergence times given in \eqref{eq:ADmC_T} along with the tracking error given in \eqref{eq:ADmC_err} for $\Pi=0$, completing the proof of the corollary. 
\end{IEEEproof}

\subsection{Exact Consensus}\label{sec:DMC:E}
\begin{Protocol}[t]
 \SetKwInOut{Input}{Input}
 \SetKwInOut{Output}{Output}
 \SetKwInOut{Set}{Set}
 \Input{Network's diameter upperbound $\Delta \in \mathbb{N}$.}
\Output{Current state $x_{i\Delta}(k)\in\mathbb{R}$ for $i\in V$.}
\FEND{
\If{$i\in A_k$\textbf{ \emph{or}} $k=0$}{
\For{$\ell =0,1,\ldots,\Delta$}{
$x^{i\ell}_{k+1}\gets u^{i\ell}_{k+1}$
}
 }
\ElseIf{$i\in R_k$}{
Gather $[x^{j0}_k,\ldots,x^{j\Delta}_k]$ from each neighbor $j \in \mathcal{N}^i$\\
Update the current states according to \\
$\displaystyle x^{i0}_{k+1} \gets u^i_k$\\
\For{$\ell =1,\ldots,\Delta$}{
$\displaystyle x^{i\ell}_{k+1} \gets \min_{j \in \mathcal{N}^i_k\bigcup \left\{i\right\}} x^{j(\ell-1)}_k$
}}}
\caption{{\small Exact Dynamic Min-Consensus (EDmC)}}
\label{alg:EDmC}
\end{Protocol}

In Protocol~\ref{alg:EDmC} we detailed the EDmC Protocol, which makes use of the following local interaction rule
\begin{align}\label{eq:protocol2}
\displaystyle x^{i0}_{k+1} &= 
\begin{cases}
u^i_{k}& \text{if } i\in R_k\\
u^i_{k+1} & \text{if } i\in A_k
\end{cases}
\\
x^{i\ell}_{k+1} &=
\begin{cases}
\displaystyle \min_{j \in \mathcal{N}^i_k\bigcup \left\{i\right\}} x^{j(\ell-1)}_k& \text{if } i\in R_k\\
u^i_{k+1} & \text{if } i\in A_k
\end{cases}\nonumber
\end{align}
with $\ell=1,\ldots,\Delta$ and $\Delta\in\mathbb{N}$ is an upperbound on the the diameter of the underlying comunication network, i.e., $\delta \geq \delta(\mathcal{G}_k)$ at any time $k\in \mathbb{N}$. At each iteration all agents gather the state values of their neighbors and update their state according to~\eqref{eq:protocol2}, which only requires local communications, while all arriving agents initialize their next state at their own inputs. Note that at the initial time $k=0$ all agents are arriving agents.

\begin{theorem}[EDmC Protocol: Tracking Error]\label{th:EDmC_track}
Consider an OMAS executing Protocol~\ref{alg:EDmC} under Assumptions~\ref{ass:slow_net}-\ref{ass:boundedu}. Consider a generic time $k_0\in\mathbb{N}$ at which the network changes, i.e., $\mathcal{G}_{k_0-1}\neq \mathcal{G}_{k_0}$.

If graph $\mathcal{G}_{k_0}$ is connected with diameter $\delta(\mathcal{G}_{k_0})$ and if 
\begin{equation}\label{eq:EDmC_delta}
\Upsilon \geq \Delta \geq \delta(\mathcal{G}_{k_0}),
\end{equation} 
Problem~\ref{prb:main} is solved with transient and convergence times
\begin{equation}\label{eq:EDmC_T}
T^c_{k_0} = T^t_{k_0} = \Delta.
\end{equation}
Moreover, if $\Upsilon \geq T^c_{k_0}$, the tracking error \eqref{eq:abserr_max} is bounded for any $k\in[k_0+T^c_{k_0},k_0+\Upsilon]$ by
\begin{equation}\label{eq:EDmC_err}
e_k \leq \varepsilon = (\Delta + 1) \Pi.
\end{equation}
\end{theorem} \medskip

\begin{IEEEproof}
The proof is in Appendix.
\end{IEEEproof}

\begin{corollary}[EDmC Protocol: Steady State Error]\label{th:EDmC_steady}
Consider an OMAS executing Protocol~\ref{alg:EDmC} under Assumption~\ref{ass:slow_net}. Consider a generic time $k_0\in\mathbb{N}$ at which the network changes, i.e., $\mathcal{G}_{k_0-1}\neq \mathcal{G}_{k_0}$ and let the inputs be constant for $k\in[k_0,k_0+\Upsilon]$.

If graph $\mathcal{G}_{k_0}$ is connected with diameter $\delta(\mathcal{G}_{k_0})$ and if 
\begin{equation*}
\Upsilon \geq \Delta \geq \delta(\mathcal{G}_{k_0}),
\end{equation*} 
Problem~\ref{prb:main} is solved with transient and convergence times
\begin{equation}\label{eq:EDMC_T_ss}
T^c_{k_0} = T^t_{k_0} = \Delta.
\end{equation}
Moreover, if $\Upsilon \geq T^c_{k_0}$, the tracking error \eqref{eq:abserr_max} is bounded for any $k\in[k_0+T^c_{k_0},k_0+\Upsilon]$ by
\begin{equation}\label{eq:EDMC_err_ss}
e_k = \varepsilon_{ss} = 0.
\end{equation}
\end{corollary} \medskip

\begin{IEEEproof} The proof is similar to the one of Corollary \ref{th:ADmC_steady}.
\end{IEEEproof}

\begin{table*}[]
 \centering
\begin{tabular}{@{}c|c|c|c|c|c@{}}
     & Transient Time & Convergence Time & Tracking Error  & Steady State Error & Conditions                                           \\
     & $T^t_{k_0}$    & $T^c_{k_0}$      & $\varepsilon$   & $\varepsilon_{ss}$ &                                                      \\ \midrule
ADMC &
  $\delta(\mathcal{G}_{k_0})$ &
  $\max \left\{ T^t_{k_0},\left\lceil \frac{\max \left\{x_{k_0}-\bar{u}_{k_0}\right\}}{\alpha - \Pi}\right\rceil\right\}$ &
  $(\delta(\mathcal{G}_{k_0})+1)\Pi+\alpha \delta(\mathcal{G}_{k_0})$ &
  $\alpha\delta(\mathcal{G}_{k_0})$ &
  $\alpha >\Pi$ \& $\Upsilon \geq \delta(\mathcal{G}_{k_0})$ \\
ADmC &
  $\delta(\mathcal{G}_{k_0})$ &
  $\max \left\{ T^t_{k_0},\left\lceil \frac{\max \left\{\ubar{u}_{k_0}-x_{k_0}\right\}}{\alpha - \Pi}\right\rceil\right\}$ &
  $(\delta(\mathcal{G}_{k_0})+1)\Pi+\alpha \delta(\mathcal{G}_{k_0})$ &
  $\alpha\delta(\mathcal{G}_{k_0})$ &
  $\alpha >\Pi$ \& $\Upsilon \geq \delta(\mathcal{G}_{k_0})$ \\
EDMC & $\Delta$       & $\Delta$         & $(\Delta+1)\Pi$ & $0$                & $\Upsilon\geq \Delta \geq \delta(\mathcal{G}_{k_0})$ \\
EDmC & $\Delta$       & $\Delta$         & $(\Delta+1)\Pi$ & $0$                & $\Upsilon\geq \Delta \geq \delta(\mathcal{G}_{k_0})$ \\ \bottomrule
\end{tabular}
\caption{Tabella riassuntiva dei risultati}
\label{tab:my-table}
\end{table*}

\section{Dynamic Network's Size Estimation}
\label{sec:applicationprotocols}

In this section we introduce two interesting problems in which the dynamic max-consensus protocols can be applied. The first problem is the one of \emph{size-estimation} of anonymous networks in which nodes can leave or join the network, thus leading to a time-varying quantity (the network's size) to be estimated. The second problem is the distributed computation of the Fiedler vector, which is the eigenvector corresponding to the smallest non-trivial eigenvalue of the graph's Laplacian matrix.

\subsection{Size-Estimation of Anonymous Network}


Here we extend the strategy proposed in \cite{Varagnolo11} for static networks to time-varying networks, i.e., nodes can leave or join the network at any time. The approach in \cite{Varagnolo11} is totally distributed and based on statistical inference concepts and can be briefly summarized as follows:
\begin{enumerate}
\item Nodes independently generate a vector of $M$ independent random numbers from a known distribution;
\item Nodes distributedly compute a specific function $f$ of all these numbers
through a consensus algorithm;
\item Each node infers the network size exploiting the statistical properties of the so computed quantities.
\end{enumerate}

In particular, we consider the max-consensus scenario, i.e., when the function $f$ of the randomly generated numbers to be estimated is the maximum of them. 
Differently from \cite{Terelius2012}, in our case nodes are allowed to leave and join the network. Thus, we extend the above scheme by adding the following rules:
\begin{itemize}
\item ADMC and EDMC Protocols are used as consensus algorithm in step $2)$;
\item When a node join the network, generate a new random number.
\end{itemize}
If the network is static (no nodes leave or join the network) the problem is the one considered with all reference signals staying constant. Since our algorithm is robust to re-initialization, every time a node leaves or joins the network, the algorithm is able to converge to the new set of inputs. Intuitively, the rate at which nodes leave or join the network is correlated to the rate of change of the maximum value of the numbers and thus one may possibly have some critical scenarios. Here, we just make the simple assumptions that our protocols can run a sufficiently high number of iterations such that an equilibrium is reached after each change of the network. We formalize this concept in the following assumption.

\begin{assumption}\label{ass:nodevarying} The minimum time $\Upsilon$ between two changes of the network ensured by Assumption~\ref{ass:slow_net} is greater or equal than the convergence time $T^c_{k_0}$ of the employed protocol, i.e.,
$$
\Upsilon \geq T^c_{k_0}
$$
\end{assumption}

Under Assumption~\ref{ass:nodevarying}, we are able to estimate and track the time varying size of the network without any synchronization among the agents, since no re-initialization is required by ADMC and EDMC Protocols. Our strategy is formalized in Protocol~\ref{alg:SEalgo}, where the following notation is used (we omit here the time-dependency $(k)$):

\begin{itemize}
\item Each agent $i$ generates $p$ random numbers $u^{i}\in[0,1]^p$ ;
\item $x^{ij}_k$ denotes the $i$-th agent's estimation of the maximum among all $u^{\ell j}_k$ for $\ell\in V$ at time $k$;
\item $\hat{n}^i_k$ is the $i$-th agent estimation of the size of the network based on estimations $x^{i\ell}_k$ for $\ell\in P$ at time $k$.  
\end{itemize}

\begin{Protocol}[t]
 \SetKwInOut{Input}{Input}
 \SetKwInOut{Output}{Output}
 \SetKwInOut{Set}{Set}
 \Input{Number of random numbers $p\in \mathbb{N}$.
 }
\Output{$\displaystyle \hat{n}^{i}_k\gets \frac{-p}{\sum_{j=1}^p\log(x^{ij}_k)}$ for $i\in V_k$.}
\FEND{
\If{$i\in A_k$\textbf{ \emph{or}} $k=0$}
{$u^{i}_{k+1}\gets $ rand$([0,1]^p)$
}
\For{$j=1,\ldots,p$}{
Update state $x^{ij}_{k+1}$ according to either Protocol~\ref{alg:ADMC} or Protocol~\ref{alg:EDMC} with inputs $[u^{1j}_k,\ldots,u^{nj}_k]$
}}
\caption{{\small Dynamic Size-Estimation (DSE)}}
\label{alg:SEalgo}
\end{Protocol}

\begin{theorem}\label{th:SEADMC}
Consider an OMAS executing Protocol~\ref{alg:SEalgo} with parameter $p\in \mathbb{N}$, $p>1$, under Assumptions~\ref{ass:slow_net}-\ref{ass:boundedu}-\ref{ass:nodevarying}. Consider a generic time $k_0\in\mathbb{N}$ at which the network changes, i.e., $\mathcal{G}_{k_0-1}\neq \mathcal{G}_{k_0}$.

If it is employed Protocol \ref{alg:ADMC} under the conditions of Corollary \ref{th:ADMC_steady}, then the expected value $\mathbb{E}[\hat{n}^i_k]$ at a steady state, i.e., $k\in[k_0+T^c_{k_0},k_0+\Upsilon]$ is 
\begin{equation}\label{eq:exvalADMC}
\mathbb{E}\left[\hat{n}^i_k\right] =  \varepsilon^{p-1}e^{\varepsilon np} (np)^{p} \Gamma(1-p,\varepsilon np),
\end{equation}

with $\varepsilon = \delta(\mathcal{G}_{k_0})\alpha$.

If it is employed Protocol \ref{alg:EDMC} under the conditions of Corollary \ref{th:EDMC_steady}, then the expected value $\mathbb{E}[\hat{n}^i_k]$ at a steady state, i.e., $k\in[k_0+T^c_{k_0},k_0+\Upsilon]$ is 
\begin{equation}\label{eq:exvalEDMC}
\mathbb{E}\left[\hat{n}^i_k\right] = \frac{np}{p-1}
\end{equation}

\end{theorem}

\begin{IEEEproof}
The proof is given in Appendix.
\end{IEEEproof}


\section{Numerical simulations} \label{sec:numsim} 

To illustrate the performance of the proposed protocol, simulation results are given in this section. First, we substantiate stability and error bounds of the proposed protocols by simulating a worst-case scenario network with line topology. Second, we applied these protocols in the context of distributed size estimation of time-varying networks, in which nodes can join and leave over time. A discussion of pros and cons of the proposed protocols is provided.

\subsection{Network with line topology}

For the sake of clarity and without loss of generality, in this subsection we limit the simulations to $[k_0,k_0+\Upsilon]$ for any $k_0\in\mathbb{N}$ in which the topology remains unchanged, according to Assumption \ref{ass:slow_net}. This allows us to show how the protocol steer the agents to track the time-varying maximum $\bar{u}_k$ value among the inputs $u^i_k$, proving the results on the transient and convergence times and on the bound on the tracking error given in Theorems \ref{th:ADMC_track}-\ref{th:EDMC_track} and Corollaries \ref{th:ADMC_steady}-\ref{th:EDMC_steady}. Dual simulations for the dynamical min-consensus problem are omitted.

\begin{figure}
\begin{center}
\includegraphics[width=0.22\textwidth]{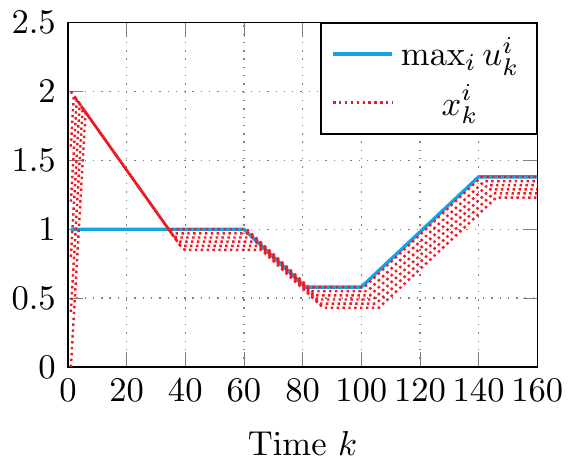}
\includegraphics[width=0.22\textwidth]{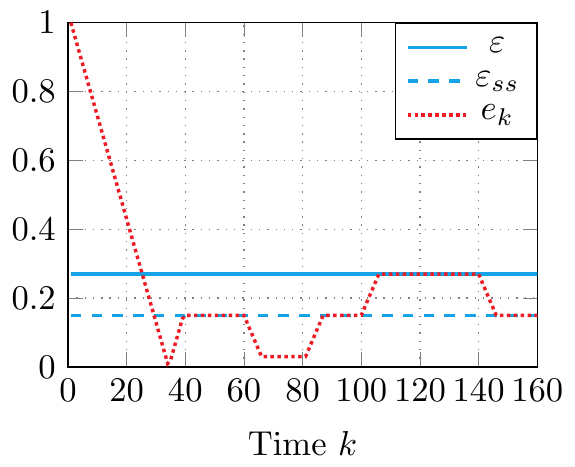}
\end{center}
\caption{Evolution of a MAS evolving according to Protocol~\ref{alg:ADMC}.}\label{fig:TwoStepPlot_P1}
\end{figure}

\begin{figure}
\begin{center}
\includegraphics[width=0.22\textwidth]{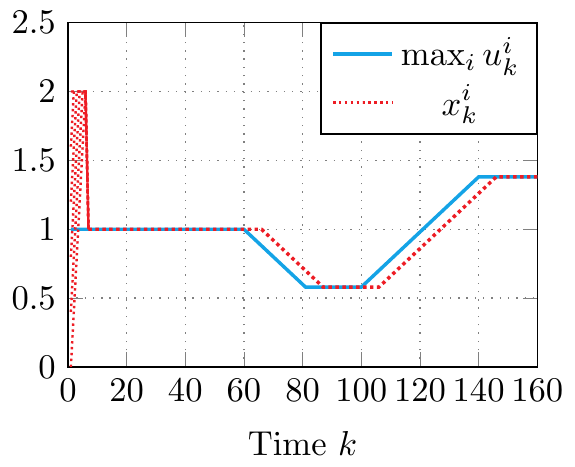}
\includegraphics[width=0.22\textwidth]{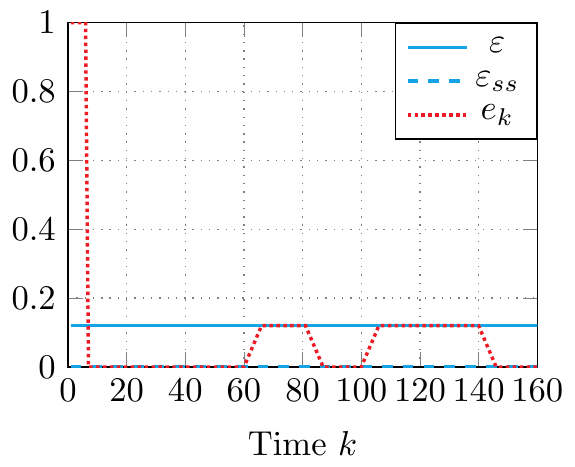}
\end{center}
\caption{Evolution of a MAS evolving according to Protocol~\ref{alg:EDMC}.}\label{fig:TwoStepPlot_P2}
\end{figure}
 
%
We simulate a network of $n=6$ agents with line topology. The choice of the line topology is instrumental to run simulations in the worst case scenario. In fact, for line graphs the information takes exactly $\delta(\mathcal{G})=n-1=5$ steps to flow through the network, thus maximizing the error for a fixed number of agents.

Figures~\ref{fig:TwoStepPlot_P1}-\ref{fig:TwoStepPlot_P2} show evolution of the state variables (dashed red lines) and maximum among the time-varying inputs (blue line) when ADMC and EDMC protocols (given in Protocol~\ref{alg:ADMC}-\ref{alg:EDMC}, respectively) are run over the MAS, respectively. State variables are initialized at $x(0)=[0,\:0.4,\:0.8,\:1.2,\:1.6,\:2]^T$ while inputs are initialized at $u(0)=[0.2,\:0.2,\:0.2,\:0.2,\:0.2,\:0.2]^T$. All inputs remain constant except for the $6$-th component, which is time-varying with respect to the following
\begin{equation}
u_6(k)
\begin{cases}
u_6(0) & \text{if }k<60\\
u_6(0)-\Pi& \text{if }k\in[60,80) \\
u_6(80)& \text{if }k\in[80,100)   \\
u_6(0)+\Pi& \text{if }k\in[100,140)\\
u_6(140) & \text{if }k\geq 140
\end{cases}\:,
\end{equation}
with $\Pi=0.02$ being the absolute change according to Assumption \ref{ass:boundedu}. 
\begin{enumerate}
\item ADMC Protocol under Theorem \ref{th:ADMC_track}:
\begin{itemize}
\item Input parameter $\alpha = 0.03$;
\item Transient time $T^t_{k_0} = 5$;
\item Convergence time $T^c_{k_0} = 34$;
\item Bound on the tracking error $\varepsilon = 0.27$;
\item Bound on the steady state error $\varepsilon_{ss} = 0.15$.
\end{itemize}
\item EDMC Protocol under Theorem \ref{th:EDMC_track}
\begin{itemize}
\item Input parameter $\Delta = 5$;
\item Transient time $T^t_{k_0} = 5$
\item Convergence time $T^c_{k_0} = 5$
\item Bound on the tracking error $\varepsilon = 0.12$;
\item Bound on the steady state error $\varepsilon_{ss} = 0$.
\end{itemize}
\end{enumerate}

\subsection{Size Estimation}


We chose to run simulations of size estimation over scale-free networks \cite{Reka2002,Amaral2000}. A scale-free network is a network whose degree distribution follows a power law, at least asymptotically. That is, the fraction $P(k)$ of nodes in the network having $k$ connections to other nodes goes for large values of $k$ as
$$
P(k)\ \sim \ k^{-\gamma}
$$
where parameter $\gamma \in \mathbb{R}$ typically is in the range $[2,3]$. Such networks are known to be \emph{ultrasmall}, as proved in \cite{Cohen2003}, meaning that their diameter scales very slow with the dimension of the network, behaving as $d \sim \ln\ln N.$

We randomly generated a scale-free network by means of Barabási–Albert (BA) model proposed in \cite{Reka2002}. This algorithm generates random scale-free networks using a preferential attachment mechanism given an initial small network, no necessarily scale-free. We use as initial network a line network of $5$ nodes, and then we run the algorithm until a network of $n=100$ nodes is generated. This network has a diameter of the order of the original small network, i.e., $d\sim 5 $. In order to simulate nodes leaving and joining the network while keeping connectivity and scale-free structure of the graph, we randomly deactivate or activate the last $m\leq 25$ nodes added to the network by the algorithm every $120$ steps. 

\begin{figure}
\begin{center}
\includegraphics[width=0.4\textwidth]{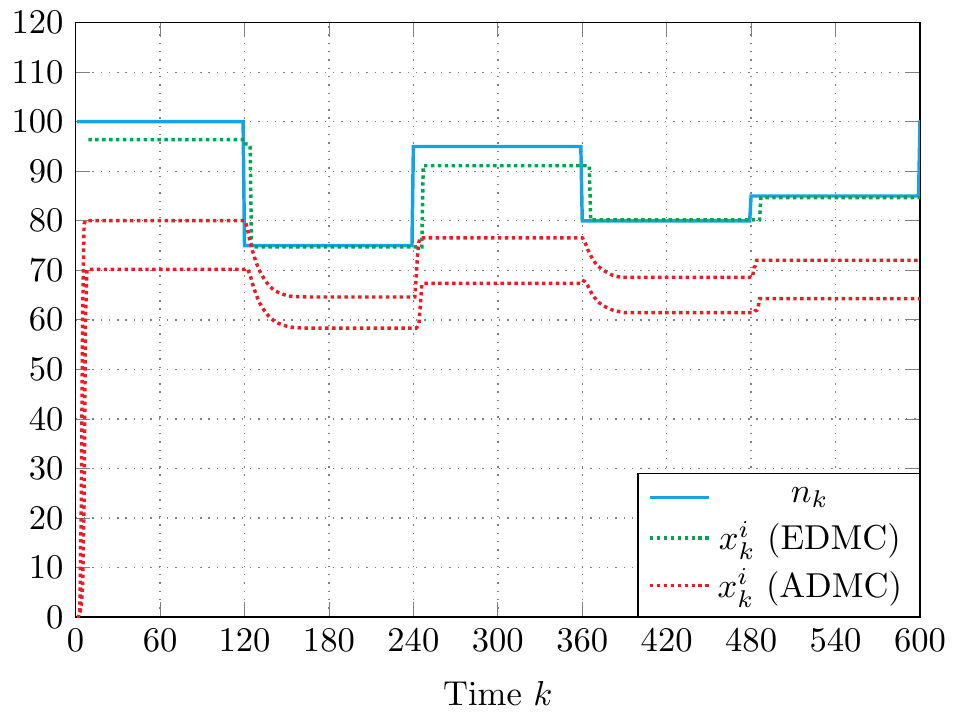}
\end{center}
\caption{Dynamic Size Estimation of a Network by means of Protocols~\ref{alg:SEalgo}.}\label{fig:SizEst}
\end{figure}

Fig.~\ref{fig:SizEst} shows the estimation of the size of a network by method proposed in \cite{Terelius2012} by means of Protocol \ref{alg:SEalgo} which makes use of one dynamic max-consensus protocols proposed in Section \ref{sec:DMC}, i.e., the ADMC protocol given in Protocol \ref{alg:ADMC} and the EDMC protocols given in Protocol \ref{alg:EDMC}.

\section{Conclusions} \label{sec:conlusions}

We have proposed, and characterized in terms time and error convergence, a distributed protocol for multi-agent systems to effectively dealing with the problem of tracking the maximum of a set of positive time-varying input reference signals. Two strengths of the proposed protocol are the following: 1) the ability to track the maximum reference signal even it is strictly lower than all states variables; 2) the robustness to initialization, meaning that the protocol is ensured to works for any initialization of the state variables. A weakness of this protocol is that exact consensus is never reached, even with constant reference inputs, thus avoiding the chance to reach a zero error. In the view of this weakness, we aim to improve the proposed protocol by means of locally distributed and time-varying tuning parameters to ensure convergence to a consensus state.

\appendix

\begin{IEEEproof}[Proof of Theorem~\ref{th:ADMC_track}]
Let us denote the maximum and, respectively, the minimum, among all agents' states, as
$$\bar{x}_k = \max_{i\in V_k} x^i_k, \qquad \ubar{x}_k = \min_{i\in V_k} x^i_k.$$
It follows that the tracking error \eqref{eq:abserr_max} satisfies 
\begin{align}
e_k & = \max_{i\in V_k} |x^i_k-\bar{u}_k|\nonumber\\
     & = \max \{|\bar{x}_k-\bar{u}_k|, |\ubar{x}_k-\bar{u}_k|\}.\label{eq:abserrgen}
\end{align}
Without loss of generality, we consider a generic time $k_0\in\mathbb{N}$ at which a change of network occurs, i.e., $\mathcal{G}_{k_0-1}\neq \mathcal{G}_{k_0}$ and $\Upsilon\geq T^c_{k_0}$. For $k\in [k_0,k_0+\Upsilon]$ the network is static. Now, for the sake of simplicity, we assume the following statements to hold, which are later proven:
\begin{equation}
\ubar{x}_k\geq \bar{u}_k-(\delta(\mathcal{G}_{k_0})+1)\Pi-\alpha\delta(\mathcal{G}_{k_0}),\quad k\in[k_0+\delta(\mathcal{G}_{k_0}),k_0+\Upsilon]\label{eq:lb_ass}.
\end{equation}
\begin{equation}
\bar{x}_k\leq\bar{u}_k+\Pi, \quad k\in[k_0+T',k_0+\Upsilon], \label{eq:ub_ass}
\end{equation}
with
\begin{equation}\label{eq:tprime}
T' = \left\lceil \frac{\max \left\{x_{k_0}-\bar{u}_{k_0}, 0\right\}}{\alpha - \Pi}\right\rceil.
\end{equation}

Relations \eqref{eq:lb_ass}-\eqref{eq:ub_ass} are both satisfied for ${T^c_{k_0} = \max\{\delta{G}_{k_0}),T'\}}$. This is the convergence time because it allows to prove the boundedness of the tracking error. In fact, for $k\in[k_0+ T^c_{k_0},k_0+\Upsilon]$ it holds
 $$x^i_k \in \left[\bar{u}_k-(\delta(\mathcal{G}_{k_0})+1)\Pi-\alpha\delta(\mathcal{G}_{k_0}),\bar{u}_k+\Pi\right].$$
and a bound on the tracking error can be computed from~\eqref{eq:abserrgen} as 
\begin{align*}
e_k & = \max_{i\in V_k} |x^i_k-\bar{u}_k|\\
     & = \max \{|\bar{x}_k-\bar{u}_k|, |\ubar{x}_k-\bar{u}_k|\}\\
     & \leq \max \{\Pi, (\delta(\mathcal{G}_{k_0})+1)\Pi+\alpha\delta(\mathcal{G}_{k_0})\}\\
     & \leq (\delta(\mathcal{G}_{k_0})+1)\Pi+\alpha \delta(\mathcal{G}_{k_0}),
\end{align*}
which corresponds to the bound given by the theorem. This proves that condition \eqref{eq:err_lay} of Problem~\ref{prb:main} is achieved.

Now, we proceed to prove that condition \eqref{eq:err_dec} of Problem~\ref{prb:main} is achieved. For $k \in [k_0+ \delta(\mathcal{G}_{k_0}),k_0+T^c_{k_0}]$, relation \eqref{eq:lb_ass} holds true but relation \eqref{eq:ub_ass} does not, thus we use the following
\begin{align}\label{eq:x+1M}
\bar{x}_{k+1}&=\max_{i\in V_k}x^i_{k+1}
\\
&=\max_{i\in V_k}\max_{j\in\mathcal{N}^i_k\cup\{i\}}\{x^j_k- \alpha,u^i_k\}
\nonumber\\
&=\max_{i\in V_k}\{x^i_k- \alpha,u^i_k\}
\nonumber\\
&=\max\{\bar{x}_k - \alpha,\bar{u}_k\},
\nonumber
\end{align}
which is due to Protocol~\ref{alg:ADMC} and its interaction rule given in eq.~\eqref{eq:protocol1}.
%
The difference in the tracking error can be computed from~\eqref{eq:abserrgen} as 
\begin{align*}
e_{k+1}-e_k & = \max_{i\in V_k} |x^i_{k+1}-\bar{u}_{k+1}| - \max_{i\in V_k} |x^i_k-\bar{u}_k|\\
     & = \max \{|\bar{x}_{k+1}-\bar{u}_{k+1}|, |\ubar{x}_{k+1}-\bar{u}_{k+1}|\}\\
     &\quad-\max \{|\bar{x}_k-\bar{u}_k|, |\ubar{x}_k-\bar{u}_k|\}.
\end{align*}
Observing that $|\ubar{x}_{k+1}-\bar{u}_{k+1}|=|\ubar{x}_k-\bar{u}_k|$ by \eqref{eq:lb_ass}, it is sufficient to prove that $|\bar{x}_{k+1}-\bar{u}_{k+1}|\leq |\bar{x}_k-\bar{u}_k|$ in order to prove $e_{k+1}-e_k\leq 0$. We compute
\begin{align*}
|\bar{x}_{k+1}-\bar{u}_{k+1}| & = |\max\{\bar{x}_k - \alpha,\bar{u}_k\}-\bar{u}_{k+1}|\\
&= \max\{|\bar{x}_k - \alpha -\bar{u}_{k+1}|,|\bar{u}_k-\bar{u}_{k+1}| \}\\
&= \max\{|\bar{x}_k - \alpha -\bar{u}_k\pm \Pi|,\Pi\}\\
&\leq  \max\{|\bar{x}_k -\bar{u}_k -\alpha|+ \Pi,\Pi\}\\
&\leq  |\bar{x}_k -\bar{u}_k -\alpha|+ \Pi
\end{align*}
For $k \in [k_0+ \delta(\mathcal{G}_{k_0}),k_0+T^c_{k_0}]$ the tracking error is surely larger than the bound $(\delta(\mathcal{G}_{k_0})+1)\Pi+\alpha \delta(\mathcal{G}_{k_0})$ and thus also than $\geq \max\{\alpha,\Pi\}$, thus one can derive
\begin{align*}
\bar{x}_k -\bar{u}_k -\alpha+ \Pi &\leq \bar{x}_k -\bar{u}_k\\
-\alpha+ \Pi &\leq 0\\
\alpha &\geq \Pi.
\end{align*}
Since by assumption it holds \eqref{eq:ADMC_alpha}, i.e., $\alpha>\Pi$ then $e_{k+1}<e_{k}$ and this proves that condition \eqref{eq:err_dec} of Problem~\ref{prb:main} is achieved and the transient time is given by
$$
T^t_{k_0} = \delta(\mathcal{G}_{k_0}).
$$
Therefore, Problem 1 is solved. To complete the proof, we proceed by proving the veracity of inequalities \eqref{eq:ub_ass}-\eqref{eq:lb_ass}.


$\bullet$ \emph{Proof of eq. \eqref{eq:lb_ass}}. It trivially holds
$$x^i_k\geq \ubar{x}_k\quad \forall i\in V_k,\:\forall k \in \mathbb{N}\:.$$
At time $k_0$ we define the set $$\mathcal{V}_0=\left\{i\in V_{k_0}: x^i_{k_0}=\bar{x}_{k_0}\right\}$$ denote the set of agents whose state at time $k_0$ is the maximum among all others. Let us now consider the set $\mathcal{V}_1$ of one-hop neighbors of nodes in set $\mathcal{V}_0$ at time $k_0+1$. Formally, $$\mathcal{V}_1=\left\{i\in V_{k_0}: (i,j)\in E, \; j \in \mathcal{V}_0 \right\}.$$
Thus, for all $i \in \mathcal{V}_1$, the state update rule~\eqref{eq:protocol1} reduces to $$x^i_{k_0+1}=\max\{\bar{x}_{k_0}-\alpha,u^i_{k_0}\},$$
because all agents $i\in \mathcal{V}_1$ have a neighbor $j\in V_0$ with state value $x^j_{k_0+1}=\bar{x}_{k_0}$. 
Thus, exploiting eq.~\eqref{eq:ub_ass} we can write
\begin{equation}
x^i_{k_0+1}\geq \bar{x}_{k_0}-\alpha, \qquad \forall i\in \mathcal{V}_1.
\end{equation}
By induction, define the set $$\mathcal{V}_{\ell}=\left\{i\in V_{k_0}: (i,j)\in E, j \in \bigcup_{s=0}^{\ell-1} \mathcal{V}_{s}\right\}.$$
It easily follows that, since the longest shortest path between two nodes in a connected graph is at most equal to its diameter $\delta(\mathcal{G}_{k_0})$, for $\ell=\delta(\mathcal{G}_{k_0})$ it holds  $\mathcal{V}_{\delta(\mathcal{G})}\equiv V_{k_0}$. Therefore
$$x^i_{k_0+\delta(\mathcal{G}_{k_0})} \geq \bar{x}_{k_0}-\delta(\mathcal{G}_{k_0})\alpha \qquad \forall i\in V_{k_0}.$$
For $k\in[k_0+\delta(\mathcal{G}_{k_0}),k_0+\Upsilon]$ we can combine the previous inequality with \eqref{eq:ub_ass} and \eqref{eq:lb_input}, leading to%
\begin{align*}
\ubar{x}(k)&\geq  \bar{x}_{k-\delta(\mathcal{G}_{k_0})} -\delta(\mathcal{G}_{k_0})\alpha\\
&\geq \bar{u}_{k-\delta(\mathcal{G}_{k_0})-1} -\delta(\mathcal{G}_{k_0})\alpha\\
&  \geq \bar{u}_k-(\delta(\mathcal{G}_{k_0})+1)\Pi-\alpha\delta(\mathcal{G}_{k_0}),
\end{align*}
which proves the veracity of eq. \eqref{eq:lb_ass}


$\bullet$ \emph{Proof of eq. \eqref{eq:ub_ass}}. Under Assumption~\ref{ass:boundedu}, at the generic time $k_0+T$ with $T\in[0,\Upsilon]$ it holds
\begin{equation}\label{eq:lb_input}
\bar{u}_k=\bar{u}_{k_0+T}\geq \bar{u}_{k_0}-T\Pi
\end{equation} and thus by \eqref{eq:x+1M} it follows
$$
\bar{x}_{k_0+T+1}=\max\{\bar{x}_{k_0} - T\alpha,\bar{u}_{k_0+T}\}.
$$
Under condition~\eqref{eq:ADMC_alpha}, one derives that the inputs vary slower than the agents' states, and therefore there exists a time $k_0+T$ after which the system \emph{reaches} the input, i.e., 
$$
\bar{x}_{k_0}-T \alpha < \bar{u}_{k_0}-T \Pi.
$$
Solving for $T$, we obtain $T'$ as in \eqref{eq:tprime}. Thus, for $k\in[k_0+T',k_0+\Upsilon]$, recalling \eqref{eq:lb_input}, the dynamics of $\bar{x}(k)$ is given by 
$$
\bar{x}_{k}= \bar{u}_{k-1} \leq \bar{u}_{k}+\Pi,
$$
proving the veracity of eq. \eqref{eq:ub_ass} and completing the proof of the theorem. 
\end{IEEEproof}

\begin{IEEEproof}[Proof of Theorem~\ref{th:EDMC_track}]
At time $k_0$, we define the set
$$
\mathcal{V}_0=\left\{i\in V_{k_0}: x^{i0}_{k_0}=\max_{j\in \mathcal{V}_{k_0}}x^{j0}_{k_0}\right\}.
$$
Since by Protocol~\ref{alg:EDMC} it holds ${x^{i0}_k = u^i_{k-1},}$ then
$$
\mathcal{V}_0=\left\{i\in \mathcal{V}_{k_0}: x^{i0}_{k_0}=\max_{j\in V_{k_0}}u^{j}_{k_0-1}\right\}.$$
Let us now consider time $k_0+1$ and the set $\mathcal{V}_1$ of one-hop neighbors of nodes in set $\mathcal{V}_0$. Formally,
$$
\mathcal{V}_1=\left\{i\in \mathcal{V}_{k_0}: (i,j)\in E, \; j \in \mathcal{V}_0 \right\}.$$
The state update rule~\eqref{eq:protocol2} for $i \in \mathcal{V}_1$, $\ell = 1$ reduces to
$$
x^{i1}_{k_0+1}=\max_{j\in\mathcal{N}^i_k\cup\{i\}}x^{j0}_{k_0}=\max_{j\in V_{k_0}}u^j_{k_0-1}=\bar{u}_{k_0-1},
$$
because all agents $i\in \mathcal{V}_1$ have a neighbor $j\in \mathcal{V}_0$ with state value $x^{j0}_{k_0}=\bar{u}_{k_0-1}$. 
By induction, for $\ell \geq 1$ define $$\mathcal{V}_{\ell}=\left\{i\in V_{k_0}: (i,j)\in E, j \in \bigcup_{s=0}^{\ell-1} \mathcal{V}_{s}\right\},$$
and therefore for all $i \in \mathcal{V}_{\ell}$ it holds 
$$
x^{i\ell}_{k_0+\ell}=\bar{u}_{k_0-1}.
$$
By noticing that $\mathcal{V}_{\Delta}=\mathcal{V}_{\delta(\mathcal{G})}\equiv V_{k_0}$, we infer that for all $i\in V_{k_0}$ and for any time $k\in[k_0 + \Delta,k_0+\Upsilon]$ with $\Upsilon \geq \Delta$, it holds \begin{equation}\label{eq:xid}
x^{i\Delta}_k=\bar{u}_{k-\Delta-1},
\end{equation}
which proves that transient and convergence times coincide and they are equal to the upperbound $\Delta$, i.e., $T^c_{k_0} = T^t_{k_0} = \Delta$. Furthermore, by Assumption~\ref{ass:boundedu} it follows 
\begin{align}\label{eq:ub}
\bar{u}_k\in[
&\bar{u}_{k-\Delta-1}-(\Delta+1)\Pi,\\\
&\bar{u}_{k-\Delta-1}+(\Delta+1)\Pi]\nonumber.
\end{align}
Finally, exploiting \eqref{eq:xid} and \eqref{eq:ub}, we conclude that for any $k\in[k_0 + \Delta,k_0+\Upsilon]$ the tracking error \eqref{eq:abserr_max} is bounded by
\begin{align*}
e_k = \max_{{i}\in V_{k_0}} |x^{i\Delta}_k - \bar{u}_k|\leq (\Delta+1)\Pi,
\end{align*}
completing the proof.

\end{IEEEproof}

\begin{IEEEproof}[Proof of Theorem~\ref{th:ADmC_track}]
We start the proof with a trivial statement
\begin{equation}\label{eq:changeA}
\ubar{u}_k= - \bar{v}_k,\qquad v_k= - u_k.
\end{equation}
Now, consider an OMAS executing Protocol~\ref{alg:ADMC} with state $y\in\mathbb{R}^n$ and inputs $v_k$, thus agents' local interaction rule is
$$
y^i_{k+1} =
\begin{cases}
\displaystyle\max_{j \in \mathcal{N}^i_k\bigcup \left\{i\right\}} \left\{ y^j_k-\alpha,v^i_k\right\}& \text{if } i\in R_k\\
v^i_{k+1} & \text{if } i\in A_k
\end{cases}
$$
Given the initial condition $y_{k_0} = - x_{k_0}$, we have that for any $k\in[k_0,k_0+\Upsilon]$ it holds for $i\in R_k$
\begin{align*}
y^i_{k+1} &= \max_{j \in \mathcal{N}^i_kj\bigcup \left\{i\right\}} \left\{ -x^j_k-\alpha,-u^i_k\right\},\\
&= - \min_{j \in \mathcal{N}^i\bigcup \left\{i\right\}} \left\{ x^j_k + \alpha,u^i_k\right\}
\end{align*}
and then by invoking \eqref{eq:changeA} it follows
$$
y^i_{k+1} = -x^i_{k+1},\quad \forall i\in V_{k_0}.
$$
It is trivial to notice that due to the above relation one can derive transient and convergence time by substitution into \eqref{eq:ADMC_T}

\begin{align*}
T^t_{k_0} &= \delta(\mathcal{G}_{k_0}),\\
T^c_{k_0} &= \max \left\{ \delta(\mathcal{G}_{k_0}),\left\lceil \frac{\max \left\{y_{k_0}-\bar{v}_{k_0}\right\}}{\alpha - \Pi}\right\rceil\right\},\nonumber\\
&= \max \left\{ \delta(\mathcal{G}_{k_0}),\left\lceil \frac{\max \left\{\ubar{u}_{k_0}-x_{k_0}\right\}}{\alpha - \Pi}\right\rceil\right\},\nonumber\\
\end{align*}
and the tracking errror by \eqref{eq:ADMC_err} as follows
\begin{align*}
e_k & = \max_{i\in V_{k_0}} |x^i_k-\ubar{u}_k|\\
	& = \max_{i\in V_{k_0}} |-y^i_k+\bar{v}_k|\\
	& = \max_{i\in V_{k_0}} |y^i_k-\bar{v}(k)|\\
     &\leq (\delta(\mathcal{G})+1)\Pi+\alpha \delta(\mathcal{G})).
\end{align*}
of Theorem~\ref{th:ADMC_track}, completing the proof.

\end{IEEEproof}

\begin{IEEEproof}[Proof of Theorem~\ref{th:EDmC_track}]
We start the proof with a trivial statement
\begin{equation}\label{eq:changeE}
\ubar{u}_k= - \bar{v}_k,\qquad v_k= - u_k.
\end{equation}
Now, consider an OMAS executing Protocol~\ref{alg:EDMC} with state $y\in\mathbb{R}^n$ and inputs $v_k$, thus agents' local interaction rule is
\begin{align*}
\displaystyle y^{i0}_{k+1} &= 
\begin{cases}
v^i_{k}& \text{if } i\in R_k\\
v^i_{k+1} & \text{if } i\in A_k
\end{cases}
\\
y^{i\ell}_{k+1} &=
\begin{cases}
\displaystyle \min_{j \in \mathcal{N}^i_k\bigcup \left\{i\right\}} y^{j(\ell-1)}_k& \text{if } i\in R_k\\
v^i_{k+1} & \text{if } i\in A_k
\end{cases}\nonumber
\end{align*}
Given the initial condition $y_{k_0} = - x_{k_0}$, we have that for any $k\in[k_0,k_0+\Upsilon]$ it holds for $i\in R_k$ and $\ell=1,\ldots,\Delta$
\begin{align*}
y^{i\ell}_{k+1} &= \max_{j \in \mathcal{N}^i\bigcup \left\{i\right\}} y^{j(\ell-1)}_k\\
&= \max_{j \in \mathcal{N}^i\bigcup \left\{i\right\}} -x^{j(\ell-1)}_k\\
&= -\min_{j \in \mathcal{N}^i\bigcup \left\{i\right\}} x^{j(\ell-1)}_k.
\end{align*}
and then by invoking \eqref{eq:changeE} it follows
$$
y^i_{k+1} = -x^i_{k+1},\quad \forall i\in V_{k_0}.
$$
It is trivial to notice that due to the above relation one can derive transient and convergence time by substitution into \eqref{eq:EDMC_T}
$$
T^t_{k_0} = T^c_{k_0} \Delta
$$
and the tracking error by \eqref{eq:EDMC_err} as follows
\begin{align*}
e_k & = \max_{i\in V_{k_0}} |x^i_k-\ubar{u}_k|\\
	& = \max_{i\in V_{k_0}} |-y^i_k+\bar{v}_k|\\
	& = \max_{i\in V_{k_0}} |y^i_k-\bar{v}(k)|\\
     &\leq (\Delta+1)\Pi.
\end{align*}
of Theorem~\ref{th:EDMC_track}, completing the proof.
\end{IEEEproof}

\begin{IEEEproof}[Proof of Theorem \ref{th:SEADMC}]
Since by Assumption~\ref{ass:slow_net} the network remains unchanged for $k\in[k_0,k_0+\Upsilon]$, then by Protocol \ref{alg:SEalgo} the inputs of the agents $V_{k_0}$ are constant in this interval. 
Thus, in the following we omit the dependence of all variables from $k$. 

For the purpose of the proof, we recall some basic concepts on order statistics \cite{}. Consider the sample $u^{1j},\ldots,u^{nj}$ consisting of the $j$-th numbers generated by the agents $i=1,\ldots,n$. The $j$-th smallest value is called the $j$-th order statistic of the sample. Let us denote $\bar{u}^j$ the $n$-th order statistics of the sample, i.e., the maximum value
$$\bar{u}_j = \max_{i \in V}u^{i j},\quad \forall j=1,\ldots,p.$$

All $u^{ij}$ are i.i.d. with probability density function $p(a)$ and and probability distribution function $P(a)$ given by
$$p(a) = 
\begin{cases}
1 & 0 \leq a \leq 1\\
0 & \text{otherwise}
\end{cases}
,\quad P(a) = 
\begin{cases}
a & 0 \leq a \leq 1\\
0 & \text{otherwise}
\end{cases},$$ 
while the probability density function of the $n$-th order statistic is given by
\begin{equation}\label{eq:pdf_nthos}
p_n(a) = nP^{n-1}(a).
\end{equation}
Consider now the sample obtained by the $n$-th order statistics of each random number generated by the agents, i.e., ${\tilde{u}=\{\bar{u}_1,\ldots,\bar{u}_p\}}$. 
Variables $\bar{u}_j$ in the sample $\tilde{u}$ depend on the parameter $n$. The likelihood function $\mathcal{L}(n|\tilde{u})$ is equal to the probability that the particular outcome $\tilde{u}$ given the parameter $n$ and, since all variables in the sample are i.i.d. random variables with probability density function \eqref{eq:pdf_nthos}, then it can be computed as the product of the probability density functions, i.e.,
$$
\mathcal{L}(n|\tilde{u}) = \prod_{j = 1}^p p_n(\bar{u}_{j}) = n^p \prod_{j = 1}^p \bar{u}_{j}^{n-1}.
$$

In practice, it is often convenient to work with the natural logarithm of the likelihood function, called the log-likelihood
\begin{align*}
\mathcal{L}^*(n|\tilde{u}) = \ln \left(\mathcal{L}(n|\tilde{u})\mathcal{}\right) &= \ln\left(n^p\prod_{j = 1}^p \bar{u}_j^{n-1}\right)\\
 &= \ln\left(n^p\right) + \ln\left(\prod_{j = 1}^p \bar{u}_j^{n-1}\right)\\
  &= \ln\left(n^p\right) + \sum_{j = 1}^p \ln\left(\bar{u}_j^{n-1}\right)\\
 &= p\ln\left(n\right) + (n-1)\sum_{j = 1}^p \ln\left(\bar{u}_j\right).
\end{align*}

The maximum likelihood estimate (MLE) is the value which maximizes log-likelihood $\mathcal{L}^*(n|x_n)$, thus giving the best estimate of $n$ from the sample $\tilde{u}$. By putting to zero the derivative in $n$, we can find an expression for the MLE,
\begin{align}\label{eq:realMLE}
MLE & = \frac{-p}{\sum_{j = 1}^p \ln\left(\bar{u}_j\right)}.
\end{align}
%
%
Equation \eqref{eq:realMLE} represents the best way to estimate the size $n$ of the network trough inference by the maximum values among the numbers generated by the agents. However, variables $\bar{u}_j$ are not known exactly at each node, and thus the best an agent can do for the estimation of $n$ is to implement the following
\begin{align}\label{eq:fakeMLE}
\widehat{MLE}_i & = \frac{-p}{\sum_{j = 1}^p \ln\left(x^{ij}\right)},\quad \forall i\in V.
\end{align}
It is necessary to understand how such the error arising from $x^{ij}\neq \bar{u}_j$ affects the estimation of $n$. We start our discussion taking into consideration the employment of Protocol \ref{alg:ADMC}.

\emph{$\bullet$ Discussion for Protocol \ref{alg:ADMC}}: By Corollary \ref{th:ADMC_steady} the steady state error in the estimating of $\bar{u}_j$ is bounded by the following 
\begin{equation}\label{eq:ejx}
e^{j} = \max_{i\in V}\left|x^{ij}-\bar{u}_j\right| \leq \:\delta(\mathcal{G})\cdot \alpha = \varepsilon.
\end{equation}
A fundamental consideration, resulting from the constructing proof of Theorem \ref{th:ADMC_track} and Corollary \ref{th:ADMC_steady}, is that at steady state the estimation $x^{ij}$ of agent $i$ of the quantity $\bar{u}_j$ is always an underestimation, i.e.,
\begin{equation*}
x^{ij}\leq \bar{u}_j,\quad \forall i\in V.
\end{equation*}
With this consideration in mind, it is easy to realize that the worst case is when at least one agent underestimates all variables $\bar{u}_j$ with maximum error \eqref{eq:ejx}. Thus, we consider such a worst case scenario by assuming that
\begin{equation}\label{eq:worstcase}
\exists i\in V:\quad  x^{ij} = \bar{u}_j-\varepsilon,\quad \forall j=1,\ldots,p.
\end{equation}
Under condition \eqref{eq:worstcase}, the worst MLE estimation $\widehat{MLE}^*$ is given by approximated MLE given by 
\begin{equation*}
\widehat{MLE}^* = \frac{-p}{\sum_{j = 1}^p \ln\left(\bar{u}_j-\varepsilon\right)}.
\end{equation*}
One can prove that condition \eqref{eq:worstcase} implies the worst case scenario since the distance $\max_{i\in V}|MLE-\widehat{MLE}_i|$ is maximized. It is straightforward to notice that $\widehat{MLE}^*\leq MLE$, thus we manipulate the above expression to get a lowerbound as follows

\begin{align}
\widehat{MLE}^* &= \frac{-p}{\sum_{j = 1}^p \ln\left(\bar{u}_j-\varepsilon\right)}\nonumber\\
&= \frac{-p}{\sum_{j = 1}^p \ln\left(\bar{u}_j(1-\frac{\varepsilon}{\bar{u}_j})\right)}\nonumber\\
&= \frac{-p}{\sum_{j = 1}^p\left[ \ln\bar{u}_j +\ln\left(1-\frac{\varepsilon}{\bar{u}_j}\right)\right]}\nonumber\\
&\geq \frac{-p}{\sum_{j = 1}^p\left[\ln\bar{u}_j +\ln\left(1-\varepsilon\right)\right]}\nonumber\\
& \geq \frac{-p}{\sum_{\ell = 1}^p \left(\ln\bar{u}_j -\varepsilon\right)}\nonumber\\
& \geq \frac{p}{\sum_{j = 1}^p (-\ln\bar{u}_j)  + p\varepsilon}\nonumber\\
& \geq \frac{1}{\frac{1}{p}\sum_{j = 1}^p (-\ln\bar{u}_j) + \varepsilon}\label{eq:lbhatn}
\end{align}

At the denominator of~\eqref{eq:lbhatn} we can recognize the term
$$\gamma =\frac{1}{p} \sum_{j=1}^p -\ln\bar{u}_j.$$
Since $(-\ln\bar{u}_j)$ with $j=1,\ldots,p$ are $p$ i.i.d. exponential random variables with rate $n$, their averaged sum $\gamma$ is known to be a gamma random variable with shape $\alpha = p$ and rate $\beta = pn$. Therefore, $\hat{n}$ is the reciprocal of a shifted gamma variable, whose probability density function $g(\cdot)$ is given by
$$
g(a)= \frac{\left(np\right)^p}{(p-1)!}a^{p-1}e^{-npa}
$$
We can now use the \emph{law of the unconscious statistician} to calculate the expected value of $\widehat{MLE}^*$ by 
\begin{equation}\label{eq:evalSEL}
\mathbb{E}[\widehat{MLE}^*]=\int_0^{\infty} f(x)g(x)dx,
\end{equation}
where $f(x) = 1/(x + \varepsilon)$.
Solution to~\eqref{eq:evalSEL} can be computed through any solver, giving the following
$$
\mathbb{E}[\widehat{MLE}^*] = \varepsilon^{p-1}e^{\varepsilon np} (np)^{p} \Gamma(1-p,\varepsilon np),
$$
where $\Gamma(p,x)$ is known as the upper incomplete gamma function. We point out that this expression holds for $n,p\in\mathbb{N}$ and $\varepsilon \in\mathbb{R}$ such that $n\geq 1$, $p>1$ and $\varepsilon\geq 0$. This completes the first part of the proof.

\emph{$\bullet$ Discussion for Protocol \ref{alg:EDMC}}: By Corollary \ref{th:EDMC_steady} the steady state error in the estimating of $\bar{u}_j$ is bounded by the following 
\begin{equation}\label{eq:ejx}
e^{j} = \max_{i\in V}\left|x^{ij}-\bar{u}_j\right| = 0 =  \varepsilon.
\end{equation}
Solution to \eqref{eq:evalSEL} for $\varepsilon=0$ can be computed by any solver, giving the following 
$$
\mathbb{E}[\widehat{MLE}^*] = \frac{np}{p-1}.
$$
We point out that this expression holds for $n,p\in\mathbb{N}$ and $\varepsilon \in\mathbb{R}$ such that $n\geq 1$, $p>1$ and $\varepsilon\geq 0$. This completes the first part of the proof. This result is coherent to the expected value with zero error given in \cite{Varagnolo11}, thus proving~\eqref{eq:exvalEDMC} and confirming that~\eqref{eq:exvalADMC} is a generalization for $\varepsilon\geq 0$.
\end{IEEEproof}

\bibliographystyle{IEEEtran}
 
\bibliography{Dynamic_Size_Estimation_arxiv}

\end{document}